\documentclass[12pt]{article}
\usepackage{amsfonts}

\usepackage{mathrsfs}
\usepackage{graphicx}
\usepackage{amsmath}

\numberwithin{equation}{section}

\def\mathcal{\mathscr}

\def\ang #1{{\langle #1\rangle}}
\newtheorem{theorem}{Theorem}[section]

\newtheorem{corollary}[theorem]{Corollary}

\newtheorem{definition}[theorem]{Definition}

\newtheorem{lemma}[theorem]{Lemma}

\newenvironment{proof}[1][Proof]{\textbf{#1.} }{\ \rule{0.5em}{0.5em}}

\begin{document}

\title{{\ {Martingale representations for diffusion processes}}\\
and backward stochastic differential equations}
\date{}
\author{\emph{{\small {\textsc{By Zhongmin Qian \ and Jiangang Ying}}}} \\
%EndAName
\emph{{\small {Oxford University and Fudan University}}}}
\maketitle

\leftskip1truecm \rightskip1truecm \noindent {\textbf{Abstract.}} In this
paper we explain that the natural filtration of a continuous Hunt process is
continuous, and show that martingales over such a filtration are continuous.
We further establish a martingale representation theorem for a class of
continuous Hunt processes under certain technical conditions. In particular
we establish the martingale representation theorem for the martingale parts
of (reflecting) symmetric diffusions in a bounded domain with a continuous
boundary. Together with an approach put forward in \cite{lyons-etc1}, our
martingale representation theorem is then applied to the study of initial
and boundary problems for quasi-linear parabolic equations by using
solutions to backward stochastic differential equations over the filtered
probability space determined by reflecting diffusions in a bounded domain
with only continuous boundary. 

\leftskip0truecm \rightskip0truecm

\tableofcontents

\newpage

\section{Introduction}

The Brownian motion, an important example of martingales, diffusions and
Gaussian processes, possesses some remarkable properties which have been the
inspiration for research in areas from probability, statistics to
mathematical finance.

Let $B=(B_{t})_{t\geq 0}$ be the Brownian motion in $\mathbb{R}^{d}$ started
with an initial distribution $\mu $. The natural filtration $(\mathscr{F}%
_{t}^{\mu })_{t\geq 0}$ (called the Brownian filtration, for a definition,
see \cite{BG1}) is continuous, and all martingales on the filtered
probability space $(\Omega ,\mathscr{F},\mathscr{F}_{t}^{\mu },P^{\mu })$
are continuous. More importantly the Brownian motion has the \emph{%
martingale representation property}: any martingale over $(\Omega ,%
\mathscr{F},\mathscr{F}_{t}^{\mu },P^{\mu })$ can be expressed as an It\^{o}
integral against Brownian motion.

The (predictable) martingale representation property of a family of
martingales has been studied in a more extended setting, and several general
results have been obtained. For example, Jacod and Yor \cite{jacod-yor1}
have discovered the equivalence between the martingale representation
property and the extremal property of martingale measures. Jacod \cite
{jacod1} also obtained the martingale representation property in terms of
the uniqueness of some martingale problems, and further present criteria in
terms of predictable characteristics. These results have greatly illuminated
the subject matter. When applied to specific situations, further work and
indeed hard estimates are often required. As a matter of fact, we still have
very limited examples of martingales and filtrations which possess
martingale representation property (see, e.g. \cite{bally}, \cite{jacod1}, 
\cite{nulart1}, \cite{zheng2}).

The renewed interest in recent years in the martingale presentation property
has been motivated not only by its own right, but also by its important
applications in the mathematical finance (\cite{Karoui-peng1}, \cite
{yong-zhou1}), backward stochastic differential equations (\cite{bally}, 
\cite{pardoux1}) and their applications in some non-linear partial
differential equations, see also for example \cite{bismut1}, \cite
{briand-hu1}, \cite{Kobylanski1}, \cite{Ma-Protter-yong1}, \cite{peng2}, 
\cite{zheng1} and the reference therein.

An intimate question is the continuity of natural filtrations generated by
semimartingales and diffusion processes. A great knowledge about them has
been obtained in the past. For example, a complete characterization of the
natural filtrations generated by simple jump processes and L\'{e}vy
processes in terms of their sample paths is known. Much information has been
obtained for a class of Markov processes. We know, from the fundamental work
by Blumenthal, Chung, Dynkin, Getoor, Hunt, Meyer etc., (see for example 
\cite{BG1}, \cite{chung1}, \cite{dynkin1}, \cite{meyer1}, in particular Hunt 
\cite{hunt1}), that the natural filtration of a Feller process is right
continuous and quasi-left continuous. On the other hand, to the best
knowledge of the present authors, there are no general conditions in
literature to guarantee the martingales over the natural filtration of a
Markov process to be continuous. A reasonable conjecture is that the natural
filtration of a diffusion process (a continuous strong Markov process)
should be continuous, so are the martingales over the natural filtration.
Such a result is plausible but remains to prove.

In this paper, we show that all martingales over the natural filtration of a 
\emph{continuous} Hunt process are continuous. Our proof follows a key idea
originated from Blumenthal \cite{Blumenthal1}, formulated carefully in Meyer 
\cite{meyer1}.

The main result of the present article is a martingale representation
theorem for a class of continuous Hunt processes which satisfy a technical
condition called \emph{the Fukushima representation property}.

As a consequence of our main result, we establish the martingale
representation theorem for symmetric diffusion processes on a domain, with
Dirichlet or Neumann boundary condition. More precisely, let $D\subset 
\mathbb{R}^{d}$ be a bounded domain with a continuous boundary $\partial D$.
Consider the symmetric diffusion in $D$ with a formal infinitesimal
generator 
\begin{equation*}
L=\frac{1}{2}\sum_{i,j=1}^{d}\frac{\partial }{\partial x^{j}}a^{ij}(x)\frac{%
\partial }{\partial x^{i}}\text{ \ \ \ in \ }D
\end{equation*}
subject to the Dirichlet or Neumann boundary condition (for precise meaning,
see \S 4 and \S 5 below), where $(a^{ij})$ are only Borel measurable and
satisfies the uniform elliptic condition. To make it clear, let 
\begin{equation*}
X=(\Omega ,\mathscr{F},\mathscr{F}_{t},X_{t},\theta _{t},P^{x})\text{.}
\end{equation*}
be the symmetric diffusion associated with the Dirichlet form $(\mathscr{E},%
\mathscr{F})$ where 
\begin{equation*}
\mathscr{E}(u,v)=\frac{1}{2}\int_{D}\sum_{i,j=1}^{d}a^{ij}\frac{\partial u}{%
\partial x^{j}}\frac{\partial v}{\partial x^{i}}dx
\end{equation*}
and the Dirichlet space $\mathscr{F}=H_{0}^{1}(D)$ or $H^{1}(D)$ depending
on the Dirichlet or Neumann boundary condition. Here we use the same letter $%
\mathscr{F}$ to denote the filtration as well as the\ Dirichlet space: we
hope it should be clear from the context which one $\mathscr{F}$ stands for. 
$\mathscr{F}$, $\mathscr{F}_{t}$ are the natural filtrations generated by
the symmetric diffusion $(X_{t})_{t\geq 0}$. It happens in this case that
the coordinate functions $u^{j}(x)=x^{j}$ belong to the local Dirichlet
space $\mathscr{F}_{\text{loc}}$, and $X_{t}=(X_{t}^{1},\cdots ,X_{t}^{d})$
has the Fukushima's decomposition 
\begin{equation*}
X_{t}^{j}-X_{0}^{j}=M_{t}^{j}+A_{t}^{j}\text{ \ \ }P^{x}\text{-a.e. }%
j=1,\cdots ,d
\end{equation*}
for all $x\in D$ (or $\overline{D}$ in the Neumann boundary condition case)
except for a zero capacity set with respect to the Dirichlet form $(%
\mathscr{E},\mathscr{F})$, where $M^{1}$, $\cdots $, $M^{d}$ are continuous
martingales additive functionals and $A^{1}$, $\cdots $, $A^{d}$ are
continuous additive functionals with zero energy. The following martingale
representation theorem follows from our main result.

\begin{theorem}
\label{mart-l1}Under the above assumptions, for any initial distribution $%
\mu $ which has no charge on capacity zero sets, the family $(M^{1}$, $%
\cdots $, $M^{d})$ of martingales over $(\Omega ,\mathscr{F}^{\mu },%
\mathscr{F}_{t}^{\mu },P^{\mu })$ has the martingale representation
property: for any square integrable martingale $N=(N_{t})_{t\geq 0}$ over $%
(\Omega ,\mathscr{F}^{\mu },\mathscr{F}_{t}^{\mu },P^{\mu })$ there are
unique $(\mathscr{F}_{t}^{\mu })$-predictable processes $F^{1}$, $\cdots $, $%
F^{d}$ such that 
\begin{equation*}
N_{t}-N_{0}=\sum_{j=1}^{d}\int_{0}^{t}F_{s}^{j}dM_{s}^{j}\text{.}
\end{equation*}
\end{theorem}

Theorem \ref{mart-l1} may be applied to the symmetric diffusions in $\mathbb{%
R}^{d}$ with Dirichlet form $(\mathscr{E},\mathscr{F})$ where $\mathscr{F}%
=H_{0}^{1}(\mathbb{R}^{d})$. This special case has been proved in \cite
{zheng2} and \cite{bally}.

In \cite{zheng2}, Zheng has pointed out that the martingale part of
symmetric diffusion $(X_{t})_{t\geq 0}$ in $\mathbb{R}^{d}$ with
infinitesimal generator being a uniform second order elliptic operator in
divergence form has the martingale predictable representation property and
described a proof based on the results on the Dirichlet process $p(t,X_{t})$
obtained in Lyons and Zheng \cite{lyons-zheng1} and \cite{lyons-zheng2},
where $p(t,x)$ is the probability density function of $X_{t}$ under the
stationary distribution. More precisely, Lyons and Zheng \cite{lyons-zheng2}
have extended Fukushima's representation theorem for martingale additive
functionals to a class of processes which has a form $f(t,X_{t})$, where $f$
has finite space-time energy. Their results in particular yield that $%
p(t,X_{t})$ is a Dirichlet process in the sense of F\"{o}llmer \cite
{follmer1}, and its martingale part can be expressed as an It\^{o} integral
against $(M^{1},\cdots ,M^{d})$, which, together with the Markov property,
allows to show that for $\xi =f_{1}(X_{t_{1}})\cdots f_{n}(X_{t_{n}})$ the
conditional expectation $E(\xi |\mathscr{F}_{t})$ is again a Dirichlet
process which can be expressed as an It\^{o}'s integral against $%
(M^{1},\cdots ,M^{d})$, where the expectation is taken against the
stationary distribution $P^{m}(\cdot )=\int_{\mathbb{R}^{d}}P^{x}(\cdot )dx$%
. A routine procedure based on the Doob's maximal inequality allows to prove
the martingale representation theorem for the symmetric diffusion in $%
\mathbb{R}^{d}$ with the generator $L$ an elliptic operator of second order.
Apparently not knowing the work \cite{zheng2}, in an independent work Bally,
Pardoux and Stoica \cite{bally}, among other things, a detailed proof has
been provided.

The technical difficulty with the proof described above lies in the fact
that even for a smooth function $f$, $f(X_{t})$ may not be a semimartingale,
so that It\^{o}'s calculus can not be applied. Instead of considering random
variables with product form such as $\xi =f_{1}(X_{t_{1}})\cdots
f_{n}(X_{t_{n}})$ which linearly span a vector space dense in $L^{p}(\Omega ,%
\mathscr{F},P^{\mu })$, we utilize a linear vector $\mathscr{C}$ spanned by
those random variables which have a product form $\xi =\xi _{1}\cdots \xi
_{n}$ with 
\begin{equation*}
\xi _{j}=\int_{0}^{\infty }e^{-\alpha _{j}t}f_{j}(X_{t})dt
\end{equation*}
where $f_{j}$ are bounded Borel measurable functions, and $\alpha _{j}>0$
for $j=1,\cdots ,n$. According to Meyer \cite{meyer1}, $\mathscr{C}$ is
dense in $L^{p}(\Omega ,\mathscr{F},P^{\mu })$. The important feature is
that $U^{\alpha }f(X_{t})$ is a semimartingale for any $\alpha >0$ and a
bounded Borel function $f$, where $U^{\alpha }$ is the resolvent of the
transition semigroup. Moreover, in the symmetric case, $U^{\alpha }f$ always
belongs to the Dirichlet space (when $f$ is bounded and square integrable)
and thus Fukushima's representation theorem for \emph{martingale additive
functionals} can be applied to extend the representation to any martingales.

The martingale representation theorem for the symmetric diffusion process $%
(X_{t})_{t\geq 0}$ in a domain $D$ allows us to study the following type of
backward stochastic differential equation (BSDE) 
\begin{equation*}
dY_{t}=-f(t,Y_{t},Z_{t})dt+\sum_{j=1}^{d}Z_{t}^{j}dM_{t}^{j}
\end{equation*}
with a terminal condition that $Y_{T}=\xi \in L^{2}(\Omega ,\mathscr{F}%
_{T}^{\mu },P^{\mu })$, and thus gives a probability representation for weak
solutions to the initial and boundary value problem for non-linear parabolic
equations. The existence and uniqueness of solutions to the BSDE follows
from exactly the same approach as for the Brownian motion case, which is the
pioneering work in BSDE done by Pardoux and Peng in \cite{pardoux1}. We
however describe an approach put forward in \cite{lyons-etc1}, which allows
us to devise an alternative probability representation for the initial and
boundary problem of the corresponding semi-linear parabolic equation 
\begin{equation*}
\frac{\partial }{\partial t}u-\frac{1}{2}\sum_{i,j=1}^{d}\frac{\partial }{%
\partial x^{j}}a^{ij}\frac{\partial }{\partial x^{i}}u+f(t,u,\nabla u)=0
\end{equation*}
subject to $\left. \frac{\partial }{\partial \nu }\right| _{\partial
D}u(t,\cdot )=0$ in a bounded domain with only continuous boundary.

The paper is organized as follows. In Section 2, we develop further an idea
put forward in Meyer \cite{meyer1}, and show both the natural filtrations
and martingales over the natural filtrations for a continuous Hunt process
are continuous. This result is hardly new but it seems not appear in the
literature yet. In order to prove this result, we devise an important while
elementary formula for $E^{\mu }(\xi |\mathscr{F}_{t}^{\mu })$, which may be
considered as a refined version of a classical formula devised firstly by
Hunt and Blumenthal for potentials and multiple potential case by Meyer. In
Section 3, we establish the main result of the paper: a martingale
representation theorem for a continuous Hunt process under technical
assumptions called the Fukushima representation property, and give some
examples in which our result may apply. In Section 4, we outline the
existence and uniqueness of solutions to backward stochastic differential
equations over the natural filtered probability space over a reflecting
symmetric diffusion in a bounded domain with non-smooth diffusion
coefficients and non-smooth boundary, and finally we apply the theory of
BSDE to the study of the initial and (Neumann) boundary problem of a
non-linear parabolic equation in a bounded domain with only continuous
boundary. We believe these results are new even for reflecting Brownian
motion in a domain with non-smooth boundary.

\section{Martingales over the filtrations of continuous Hunt processes}

Consider a Markov process $(\Omega ,\mathscr{F}^{0},\mathscr{F}%
_{t}^{0},X_{t},\theta _{t},P^{x})$ in a state space $E^{\prime }=E\cup
\{\partial \}$, where $E$ is a locally compact separable metric space $E$,
with transition probability function $\{P_{t}(x,\cdot ):t\geq 0\}$, i.e., 
\begin{equation}
E^{x}\{f(X_{t+s})|\mathscr{F}_{s}^{0}\}=\int_{E}f(z)P_{t}(X_{s},dz)\text{ .}
\label{es-j1}
\end{equation}
In (\ref{es-j1}), $E^x$ on the left-hand side stands for the (conditional)
expectation with respect to the probability measure $P^{x}$, and the
right-hand side may be abbreviated as $P_{t}f(X_{s})$ where 
\begin{equation*}
P_{t}f(x)=\int_{E}f(z)P_{t}(x,dz)
\end{equation*}
which is well defined for a bounded or non-negative Borel measurable
function $f$. The family of kernels $(P_{t})_{t>0}$ is called the transition
semigroup associated with the Markov process $(X_{t})_{t\geq 0}$. Without
specification, $(\mathscr{F}_{t}^{0})_{t\geq 0}$ is the filtration generated
by $(X_{t})_{t\geq 0}$, that is $\mathscr{F}_{t}^{0}=\sigma \{X_{s}:s\leq
t\} $ \ for $t\geq 0$ and $\mathscr{F}^{0}=\sigma \{X_{s}:s\geq 0\}$.

For a $\sigma $-finite measure $\mu $ on $(E,\mathscr{B}(E))$ (where $%
\mathscr{B}(M)$ always represents the Borel $\sigma $-algebra on a
topological space $M$) 
\begin{equation*}
P^{\mu }(\Lambda)=\int_{E}P^{x}(\Lambda)\mu (dx),\ \ \Lambda\in\mathscr{F}^0,
\end{equation*}
defines a measure on $(\Omega ,\mathscr{F}^{0})$. If $\mu $ is a
probability, then $(X_{t})_{t\geq 0}$ is Markovian under $P^{\mu }$ with
transition semigroup $(P_{t})_{t>0}$ and initial distribution $\mu $, in the
sense that 
\begin{equation*}
E^{\mu }\{f(X_{t+s})|\mathscr{F}_{s}^{0}\}=\int_{E}f(z)P_{t}(X_{s},dz)\text{
\ \ }P^{\mu }\text{-a.e.}
\end{equation*}
and $P^{\mu }\{X_{0}\in A\}=\mu (A)$ for any $A\in \mathscr{B}(E)$, where $%
E^{\mu }$ is the (conditional) expectation against the probability $P^{\mu }$%
.

Denote by $\mathscr{P}(E)$ the space of all probability measures on $(E,%
\mathscr{B}(E))$. If $\mu \in \mathscr{P}(E)$, $\mathscr{F}^{\mu }$ denotes
the completion of $\mathscr{F}^{0}$ under $P^{\mu }$, and $\mathscr{F}%
_{t}^{\mu }$ is the smallest $\sigma $-algebra containing $\mathscr{F}%
_{t}^{0}$ and all sets in $\mathscr{F}^{\mu }$ with zero probability. $(%
\mathscr{F}_{t}^{\mu })_{t\geq 0}$ is called the natural filtration of the
Markov process $(X_{t})_{t\geq 0}$ with initial distribution $\mu $. Let 
\begin{equation*}
\mathscr{F}_{t}=\bigcap_{\mu \in \mathscr{P}(E)}\mathscr{F}_{t}^{\mu }
\end{equation*}
and $(\mathscr{F}_{t})_{t\geq 0}$ is called the natural filtration
determined by the Markov process 
\begin{equation*}
X=(\Omega ,\mathscr{F}^{0},\mathscr{F}_{t}^{0},X_{t},\theta _{t},P^{x}).
\end{equation*}

If $(\mathscr{G}_{t})_{t\geq 0}$ is a filtration, i.e. an increasing family
of $\sigma $-algebras on a common sample space, then $\mathscr{G}_{t+}=\cap
_{s>t}\mathscr{G}_{s}$ for $t\geq 0$ and $\mathscr{G}_{t-}=\sigma \{%
\mathscr{G}_{s}:s<t\}$ for $t>0$. The filtration is called right (resp.
left) continuous if $\mathscr{G}_{t+}=\mathscr{G}_{t}$ for all $t\geq 0$
(resp. $\mathscr{G}_{t-}=\mathscr{G}_{t}$ for all\thinspace $t>0$). The
sample function properties of a Markov process $(X_{t})_{t\geq 0}$ and the
continuity properties of its natural filtration had been studied by
Blumenthal, Dynkin, Getoor, Hunt, Meyer etc. The fundamental results have
been established via the regularity of the transition probability function $%
\{P_{t}(x,\cdot ):t>0\}$. Their work achieved the climax for Markov
processes with Feller transition semigroups.

As matter of fact, the continuity of the filtration $(\mathscr{F}_{t}^{\mu
}) $ (or $(\mathscr{F}_{t})$) does not follow that of sample function $%
(X_{t})_{t\geq 0}$. For example, a right continuous Markov process does not
necessarily lead to the right continuity of its natural filtration $(%
\mathscr{F}_{t}^{\mu })$ (or $(\mathscr{F}_{t})$). The same claim applies to
the left continuity. In fact, the regularity of natural filtrations is much
to do with the nature of the Markov property, such as strong Markov property.

Let $C_{\infty }(E)$ (resp. $C_{0}(E)$) denote the space of all continuous
functions $f$ on $E$ which vanish at infinity $\partial $, i.e. $%
\lim_{x\rightarrow \partial }f(x)=0$ (resp. with compact support). Recall
that a transition semigroup $(P_{t})_{t>0}$ on $(E,\mathscr{B}(E))$ is
Feller, if for each $t>0$, $P_{t}$ preserves $C_{\infty }(E)$ and $%
\lim_{t\downarrow 0}P_{t}f(x)=f(x)$ for each $x\in E$ and $f\in $ $C_{\infty
}(E)$.

For a given Feller semigroup $(P_{t})_{t>0}$ on $(E,\mathscr{B}(E))$, there
is a Markov process 
\begin{equation*}
X=(\Omega ,\mathscr{F}^{0},\mathscr{F}_{t}^{0},X_{t},\theta _{t},P^{x})
\end{equation*}
with the\ Feller transition semigroup $(P_{t})_{t>0}$ such that the sample
function $t\rightarrow X_{t}$ is right continuous on $[0,\infty )$ with left
hand limits on $(0,\infty )$. In this case, we call $(\Omega ,\mathscr{F}%
^{0},\mathscr{F}_{t}^{0},X_{t},\theta _{t},P^{x})$ a Feller process on $E$.

For a Feller process $(\Omega ,\mathscr{F}^{0},\mathscr{F}%
_{t}^{0},X_{t},\theta _{t},P^{x})$, the natural filtration $(\mathscr{F}%
_{t}^{\mu })_{t\geq 0}$ for any $\mu \in \mathscr{P}(E)$ and as well as $(%
\mathscr{F}_{t})_{t\geq 0}$ are right continuous. $\ (X_{t})_{t\geq 0}$ and $%
(\mathscr{F}_{t}^{\mu })_{t\geq 0}$ are also quasi-left continuous, that is,
if $T_{n}$ is an increasing family of $(\mathscr{F}_{t}^{\mu })$-stopping
times, and $T_{n}\uparrow T$, then $\lim_{n\rightarrow \infty
}X_{T_{n}}=X_{T}$ on $\{T<\infty \}$ and $\mathscr{F}_{T}^{\mu }=\sigma \{%
\mathscr{F}_{T_{n}}^{\mu }:n\in \mathbb{N}\}$. Therefore accessible $(%
\mathscr{F}_{t}^{\mu })$-stopping times are predictable. An $(\mathscr{F}%
_{t}^{\mu })$-stopping time $T$ is totally inaccessible if and only if $%
P^{\mu }\{T<\infty \}>0$ and $X_{T}\neq X_{T-}$ on $\{T<\infty \}$ $P^{\mu }$%
-a.e. Similarly, $T$ is accessible if and only if $X_{T}=X_{T-}$ $P^{\mu }$%
-a.e. on $\{T<\infty \}$. Hence $X$ has only inaccessible jump times.

What we are mainly concerned in this article is Hunt processes. Hunt
processes are right continuous, strong Markov processes which are quasi-left
continuous. These processes are defined in terms of sample functions, rather
than transition semigroups, see \cite{BG1} and \cite{chung1} for details. It
is well-known that Feller processes are stereotype of Hunt processes or the
later is an abstraction of the former.

We are interested in the martingales over the filtered probability space $%
(\Omega ,\mathscr{F}^{\mu },\mathscr{F}_{t}^{\mu },P^{x})$, and we are going
to show that, if $(\Omega ,\mathscr{F}^{0},\mathscr{F}_{t}^{0},X_{t},\theta
_{t},P^{x})$ is a Hunt process which has continuous sample function, then
any martingale on this filtered probability space is continuous, a result
one could expect for the natural filtration of a diffusion process. Indeed,
this result was proved more or less by Meyer in his Lecture Notes in
Mathematics 26, ``Processus de Markov''. Meyer himself credited his proof to
Blumenthal and Getoor, more precisely a calculation done by Blumenthal \cite
{Blumenthal1}. However it is surprising that the full computation, which
yields more information about martingales over the natural filtration of a
Hunt process, was not reproduced either in the new edition of Meyer's
``Probabilit\'{e}s et Potentiels'' or Chung's ``Lectures from Markov
Processes to Brownian Motion'', although it was mentioned in \cite
{chung-walsh1} where Chung and Walsh gave an alternative proof of Meyer's
predictability result, so that Blumenthal's computation is no longer needed.
However, it is fortunate that Blumenthal's calculation indeed leads to a
proof of a martingale representation theorem we are going to establish for
certain Hunt processes, see section \S 3.

Let us first describe an elementary calculation, originally according to
Meyer \cite{meyer1} due to Blumenthal. Let 
\begin{equation*}
X=(\Omega ,\mathscr{F}^{0},\mathscr{F}_{t}^{0},X_{t},\theta _{t},P^{x})
\end{equation*}
be a Hunt process in a state space $E^{\prime }=E\cup \{\partial \}$ with
the transition semigroup $(P_{t})_{t\geq 0}$, where $\partial $ plays a role
of cemetery. Let $\{U^{\alpha }:\alpha >0\}$ be the resolvent of the
transition semigroup $(P_{t})_{t\geq 0}$: 
\begin{equation*}
U^{\alpha }(x,A)=\int_{0}^{\infty }e^{-\alpha t}P_{t}(x,A)dt
\end{equation*}
and $(U^{\alpha })_{\alpha >0}$ the corresponding resolvent (operators),
i.e. 
\begin{eqnarray*}
U^{\alpha }f(x) &=&\int_{E}f(z)U^{\alpha }(x,dz) \\
&=&\int_{0}^{\infty }e^{-\alpha t}P_{t}f(x)dt
\end{eqnarray*}
for bounded or nonnegative Borel measurable function $f$ on $E$. To save
words, we use $\mathscr{B}_{b}(E)$ to denote the algebra of all bounded
Borel measurable functions on $E$. Obviously, $C_{\infty }(E)\subset %
\mathscr{B}_{b}(E)$.

Let $K(E)\subset \mathscr{B}_{b}(E)$ be a vector space which generates the
Borel $\sigma $-algebra $\mathscr{B}(E)$. Let $\mathscr{C}\subset
L^{1}(\Omega ,\mathscr{F}^{\mu },P^{\mu })$ (for any initial distribution $%
\mu $) be the vector space spanned by all $\xi =\xi _{1}\cdots \xi _{n}$ for
some $n\in \mathbb{N}$, 
\begin{equation*}
\xi _{j}=\int_{0}^{\infty }e^{-\alpha _{j}t}f_{j}(X_{t})dt\text{ }
\end{equation*}
where $\alpha _{j}$ are positive numbers, $f_{j}\in K(E)$, $j=1,\cdots ,n$.
Meyer \cite{meyer1} proved that $\mathscr{C}$ is dense in $L^{1}(\Omega ,%
\mathscr{F}^{\mu },P^{\mu })$ for a Hunt process. Since this density result
will play a crucial role in what follows, we include Meyer's a proof for
completeness and for the convenience of the reader.\ The key observation in
the proof is the following result from real analysis.

\begin{lemma}
Let $T>0$. Let $\mathbb{K}$ denote the vector space spanned by all functions 
$e_{\alpha }(t)=e^{-\alpha t}$, where $\alpha >0$, then $\mathbb{K}$ is
dense in $C[0,T]$ equipped with the uniform norm.
\end{lemma}

The lemma follows from \ Stone-Weierstrass' theorem.

\begin{lemma}[P. A. Meyer]
For any initial distribution $\mu $ and $p\in \lbrack 1,\infty )$, $%
\mathscr{C}$ is dense in $L^{p}(\Omega ,\mathscr{F}^{\mu },P^{\mu })$.
\end{lemma}

\begin{proof}
First, by utilizing Doob's martingale convergence theorem, it is easy to
show that the collection $\mathscr{A}$ of all random variables which have
the following form 
\begin{equation*}
g_{1}(X_{t_{1}})\cdots g_{n}(X_{t_{n}})\text{,}
\end{equation*}
where $n\in \mathbb{N}$, $0<t_{1}<\cdots <t_{n}<\infty $ and $g_{j}\in K(E)$%
, is dense in $L^{p}(\Omega ,\mathscr{F}^{\mu },P^{\mu })$. Let $\mathscr{H}$
be the linear space spanned by all $\xi =\eta _{1}\cdots \eta _{n}$, where 
\begin{equation*}
\text{ }\eta _{j}=\int_{0}^{\infty }g_{j}(X_{t})\varphi _{j}(t)dt\text{,}
\end{equation*}
where $g_{j}\in K(E)$ and $\varphi _{j}\in C[0,\infty )$ with compact
supports. According to the previous lemma, for every $\varepsilon >0$ we may
choose $\psi _{j}\in \mathbb{K}$ such that 
\begin{equation*}
|\varphi _{j}(t)-\psi _{j}(t)|<\varepsilon e^{-\lambda t}\text{ \ \ for all }%
t\geq 0
\end{equation*}
for some $\lambda >0$. Let 
\begin{equation*}
\xi _{j}=\int_{0}^{\infty }g_{j}(X_{t})\psi _{j}(t)dt\text{.}
\end{equation*}
Then $\tilde{\xi}=\xi _{1}\cdots \xi _{n}\in \mathscr{C}$, and 
\begin{equation*}
|\eta _{j}-\xi _{j}|\leq \frac{1}{\lambda }||g_{j}||_{\infty }\varepsilon
\end{equation*}
where $||\cdot ||_{\infty }$ is the supermum norm. It follows that 
\begin{equation*}
E|\xi -\tilde{\xi}|^{p}\leq \frac{n^{p}}{\lambda ^{p}}\max_{j}||g_{j}||_{%
\infty }^{np}\varepsilon ^{p}
\end{equation*}
and thus $\xi $ belongs to the closure of $\mathscr{C}$. Finally it is clear
that any element 
\begin{equation*}
g_{1}(X_{t_{1}})\cdots g_{n}(X_{t_{n}})=\lim_{k\rightarrow \infty }\eta
_{1}^{k}\cdots \eta _{n}^{k}
\end{equation*}
where 
\begin{equation*}
\eta _{j}^{k}=\int_{0}^{\infty }g_{j}(X_{t})\varphi _{j}^{k}(t)dt
\end{equation*}
and $\varphi _{j}^{k}$ has compact support and $\varphi _{j}^{k}\rightarrow
\delta _{t_{j}}$ weakly. We thus have completed the proof.
\end{proof}

\medskip

Let $\mu $ be any fixed initial distribution. If $f$ is a bounded Borel
measurable function on $E$ and $\alpha >0$, then 
\begin{equation*}
\xi =\int_{0}^{\infty }e^{-\alpha t}f(X_{t})dt\in L^{1}(\Omega ,\mathscr{F}%
^{\mu },P^{\mu })\text{.}
\end{equation*}
Consider the martingale $M_{t}=E^{\mu }\left\{ \xi |\mathscr{F}_{t}^{\mu
}\right\} $ where $t\geq 0$. According to an elementary computation in the
theory of Markov processes, 
\begin{eqnarray*}
M_{t} &=&E^{\mu }\left\{ \int_{0}^{\infty }e^{-\alpha s}f(X_{s})ds|%
\mathscr{F}_{t}^{\mu }\right\}  \\
&=&\int_{0}^{t}e^{-\alpha s}f(X_{s})ds+E^{\mu }\left\{ \int_{t}^{\infty
}e^{-\alpha s}f(X_{s})ds|\mathscr{F}_{t}^{\mu }\right\}  \\
&=&\int_{0}^{t}e^{-\alpha s}f(X_{s})ds+e^{-\alpha t}\int_{0}^{\infty
}e^{-\alpha s}P_{s}f(X_{t})ds \\
&=&\int_{0}^{t}e^{-\alpha s}f(X_{s})ds+e^{-\alpha t}U^{\alpha }f(X_{t})\text{%
.}
\end{eqnarray*}
It is known that if $X=(X_{t})_{t\geq 0}$ is a Hunt process,%
% with \emph{%continuous} sample function,
then for any $\alpha >0$ and bounded Borel measurable function $f$, $%
U^{\alpha }f$ is finely continuous, i.e., 
%(more precisely has a fine continuous modification,
$t\rightarrow U^{\alpha }f(X_{t})$ is right continuous. Moreover if $X$ is a
continuous Hunt process, 
% (or more precisely has a continuous modification), and,
it follows from a result proved by Meyer that $t\rightarrow U^{\alpha
}f(X_{t})$ is continuous, and therefore, the martingale $M_{t}=E^{\mu }\{\xi
|\mathscr{F}_{t}^{\mu }\}$ is continuous. 
%Meyer's result which seems not well known, so that
We record Meyer's result as a lemma here. This result was proved in \cite
{meyer1} for Hunt processes (see T15 THEOREME, page 89, \cite{meyer1}). A
simpler proof for\ Feller processes may be found on page 168, \cite
{Dellacherie-M1}.

\begin{lemma}[Meyer]
\label{meyer-lemma} Let $(X_{t})_{t\geq 0}$ be a Hunt process, $f\in %
\mathscr{B}_{b}(E)$, $\alpha >0$ and $h=U^{\alpha }f$ be a potential. Then 
\begin{equation*}
h(X_{t-})=h(X)_{t-}\text{ \ \ \ \ \ }\forall t>0\text{ \ \ }P^{\mu }\text{%
-a.e.}
\end{equation*}
for any initial distribution $\mu $.
\end{lemma}

P. A. Meyer pointed out that the previous computation can be carried out
equally for random variables on $(\Omega ,\mathscr{F}^{0})$ which have a
product form $\xi =\xi _{1}\cdots \xi _{n}$ where each $\xi
_{j}=\int_{0}^{\infty }e^{-\alpha _{j}s}f_{j}(X_{s})ds$. Let $\pi _{n}$
denote the permutation group of $\{1,\cdots ,n\}$.

\begin{lemma}[Blumenthal and Meyer]
\label{l:BM} Let $\xi =\xi _{1}\cdots \xi _{n}$ where $\xi
_{j}=\int_{0}^{\infty }e^{-\alpha _{j}s}f_{j}(X_{s})ds$, $\alpha _{j}>0$ and 
$f_{j}\in \mathscr{B}_{b}(E)$, and $M_{t}=E^{\mu }\left\{ \xi |\mathscr{F}%
_{t}^{\mu }\right\} $. Then 
\begin{align}
M_{t} =\sum_{k=0}^{n}\sum_{(j_{1},\cdots ,j_{k},\cdots ,j_{n})\in \pi
_{n}}\left( \prod_{i=1}^{k}\int_{0}^{t}e^{-\alpha
_{j_{i}}s}f_{j_{i}}(X_{s})ds\right) \cdot F_{(j_{1},\cdots ,j_{k},\cdots
,j_{n})}(X_{t})  \label{for-01}
\end{align}
where 
\begin{equation}
F_{(j_{1},\cdots ,j_{k},\cdots ,j_{n})}(x)=E^{x}\left\{ \left(
\prod_{l=k+1}^{n}e^{-\alpha _{j_{l}}t}\int_{0}^{\infty }e^{-\alpha
_{j_{l}}s}f_{j_{l}}(X_{s})ds\right) \right\} \text{.}  \label{for-02}
\end{equation}
\end{lemma}

\begin{proof}
The task is to calculate the conditional expectation $M_{t}=E^{\mu }\left\{
\xi |\mathscr{F}_{t}^{\mu }\right\} $. The idea is very simple: spliting
each $\xi _{j}$ into 
\begin{equation*}
\xi _{j}=\int_{0}^{t}e^{-\alpha _{j}s}f_{j}(X_{s})ds+e^{-\alpha
_{j}t}\int_{0}^{\infty }e^{-\alpha _{j}s}f_{j}(X_{s}\circ \theta _{t})ds
\end{equation*}
so that 
\begin{eqnarray*}
\xi &=&\sum_{k=0}^{n}\sum_{(j_{1},\cdots ,j_{k},\cdots ,j_{n})\in \pi
\{1,\cdots ,n\}}\left( \prod_{i=1}^{k}\int_{0}^{t}e^{-\alpha
_{j_{i}}s}f_{j_{i}}(X_{s})ds\right) \\
&&\cdot \left( \prod_{l=k+1}^{n}e^{-\alpha _{j_{l}}t}\int_{0}^{\infty
}e^{-\alpha _{j_{l}}s}f_{j_{l}}(X_{s}\circ \theta _{t})ds\right) \text{.}
\end{eqnarray*}
By using the Markov property one thus obtains 
\begin{eqnarray}
M_{t} &=&E^{\mu }\left\{ \xi |\mathscr{F}_{t}^{\mu }\right\}  \notag \\
&=&\sum_{k=1}^{n}\sum_{(j_{1},\cdots ,j_{k},\cdots ,j_{n})\in \pi \{1,\cdots
,n\}}\left( \prod_{i=1}^{k}\int_{0}^{t}e^{-\alpha
_{j_{i}}s}f_{j_{i}}(X_{s})ds\right)  \notag \\
&\ \ &\cdot E^{\mu }\left\{ \left( \prod_{l=k+1}^{n}e^{-\alpha
_{j_{l}}t}\int_{0}^{\infty }e^{-\alpha _{j_{l}}s}f_{j_{l}}(X_{s}\circ \theta
_{t})ds\right) |\mathscr{F}_{t}^{\mu }\right\}  \notag \\
&=&\sum_{k=0}^{n}\sum_{(j_{1},\cdots ,j_{k},\cdots ,j_{n})\in \pi \{1,\cdots
,n\}}\left( \prod_{i=1}^{k}\int_{0}^{t}e^{-\alpha
_{j_{i}}s}f_{j_{i}}(X_{s})ds\right)  \notag \\
&\ \ &\cdot E^{X_{t}}\left\{ \left( \prod_{l=k+1}^{n}e^{-\alpha
_{j_{l}}t}\int_{0}^{\infty }e^{-\alpha _{j_{l}}s}f_{j_{l}}(X_{s})ds\right)
\right\} \text{.}  \label{jk-01}
\end{eqnarray}
\end{proof}

Our only contribution in this aspect is the following formula, which allows
to prove not only that all martingales over the natural filtration of a
continuous Hunt process are continuous, but also a martingale representation
theorem in the next section.

\begin{lemma}
\label{lem1}Let $\alpha _{j}$ be positive numbers and $f_{j}\in $ $%
\mathscr{B}_{b}(E)$ for $j=1,\cdots ,k$. Consider 
\begin{eqnarray*}
F(x) &=&\idotsint_{0<s_{1}<\cdots <s_{k}<\infty }e^{-\sum_{j=1}^{k}\alpha
_{j}s_{j}}\int_{E^{\otimes k}}f_{1}(z_{1})\cdots
f_{k}(z_{k})P_{s_{1}}(x,dz_{1}) \\
&\ &\cdot P_{s_{2}-s_{1}}(z_{1},dz_{2})\cdots
P_{s_{k}-s_{k-1}}(z_{k-1},dz_{k})ds_{1}\cdots ds_{k}\text{.}
\end{eqnarray*}
Then 
\begin{equation}
F=U^{\alpha _{1}+\cdots +\alpha _{k}}\left( f_{1}(U^{\alpha _{2}+\cdots
+\alpha _{k}}f_{2}\cdots (U^{\alpha _{k}}f_{k})\cdots \right) \text{.}
\label{re-01}
\end{equation}
\end{lemma}

\begin{proof}
To see why it is true, we begin with the case that $k=1$. In this case $%
F=\int_{0}^{\infty }e^{-\alpha s}P_{s}fds=U^{\alpha }f$. If $k=2$, then 
\begin{eqnarray*}
F &=&\iint_{0<s_{1}<s_{2}<\infty }e^{-\alpha _{2}s_{2}}e^{-\alpha
_{1}s_{1}}P_{s_{1}}\left( f_{1}P_{s_{2}-s_{1}}f_{2}\right) ds_{1}ds_{2} \\
&=&\int_{0}^{\infty }\int_{t}^{\infty }e^{-\alpha _{1}t}e^{-\alpha
_{2}s}P_{t}\left( f_{1}P_{s-t}f_{2}\right) dsdt \\
&=&\int_{0}^{\infty }e^{-\alpha _{1}t}P_{t}\left( \int_{t}^{\infty
}e^{-\alpha _{2}s}f_{1}P_{s-t}f_{2}ds\right) dt \\
&=&\int_{0}^{\infty }e^{-\alpha _{1}t}e^{-\alpha _{2}t}P_{t}\left(
f_{1}\int_{0}^{\infty }e^{-\alpha _{2}s}P_{s}f_{2}ds\right) dt \\
&=&\int_{0}^{\infty }e^{-(\alpha _{1}+\alpha _{2})t}P_{t}\left(
f_{1}U^{\alpha _{2}}f_{2}\right) dt \\
&=&U^{\alpha _{1}+\alpha _{2}}\left( f_{1}U^{\alpha _{2}}f_{2}\right)
\end{eqnarray*}
and by an induction argument, for a general case. Indeed, if $k>2$, then 
\begin{eqnarray*}
F(x) &=&\idotsint_{0<s_{1}<\cdots <s_{k+1}<\infty }e^{-\sum_{j=1}^{k}\alpha
_{j}s_{j}}e^{-\alpha _{k+1}s_{k+1}} \\
&&\times \int_{E^{\otimes k}}f_{1}(z_{1})\cdots \left(
f_{k}(z_{k})P_{s_{k+1}-s_{k}}f_{k+1}(z_{k})\right) \\
&&\times P_{s_{1}}(x,dz_{1})P_{s_{2}-s_{1}}(z_{1},dz_{2})\cdots
P_{s_{k}-s_{k-1}}(,z_{k-1},dz_{k})ds_{1}\cdots ds_{k}ds_{k+1} \\
&=&\idotsint_{0<s_{1}<\cdots <s_{k}<\infty }e^{-\sum_{j=1}^{k}\alpha
_{j}s_{j}}\int_{E^{\otimes k}}f_{1}(z_{1})\cdots \\
&&\times \left( f_{k}(z_{k})\int_{s_{k}}^{\infty }e^{-\alpha
_{k+1}s_{k+1}}P_{s_{k+1}-s_{k}}f_{k+1}(z_{k})ds_{k+1}\right) \\
&&\times P_{s_{1}}(x,dz_{1})P_{s_{2}-s_{1}}(z_{1},dz_{2})\cdots
P_{s_{k}-s_{k-1}}(,z_{k-1},dz_{k})ds_{1}\cdots ds_{k} \\
&=&\idotsint_{0<s_{1}<\cdots <s_{k}<\infty }e^{-\sum_{j=1}^{k}\alpha
_{j}s_{j}}\int_{E^{\otimes k}}f_{1}(z_{1})\cdots \\
&&\times \left( f_{k}(z_{k})e^{-\alpha _{k+1}s_{k}}\int_{0}^{\infty
}e^{-\alpha _{k+1}t}P_{t}f_{k+1}(z_{k})dt\right) \\
&&\times P_{s_{1}}(x,dz_{1})P_{s_{2}-s_{1}}(z_{1},dz_{2})\cdots
P_{s_{k}-s_{k-1}}(,z_{k-1},dz_{k})ds_{1}\cdots ds_{k} \\
&=&\idotsint_{0<s_{1}<\cdots <s_{k}<\infty }e^{-\sum_{j=1}^{k-1}\alpha
_{j}s_{j}-(\alpha _{k+1}+\alpha _{k})s_{k}} \\
&&\times \int_{E^{\otimes k}}f_{1}(z_{1})\cdots f_{k}(z_{k})U^{\alpha
_{k+1}}f_{k+1}(z_{k})dt
\end{eqnarray*}
and the formula follows the induction assumption.
\end{proof}

\begin{lemma}
\label{lem2}Let $f_{1},\cdots ,f_{k}\in \mathscr{B}_{b}(E)$, $\alpha _{j}$
positive numbers, and 
\begin{equation*}
F(x)=E^{x}\left\{ \left( \prod_{j=1}^{k}\int_{0}^{\infty }e^{-\alpha
_{j}s}f_{j}(X_{s})ds\right) \right\} \text{.}
\end{equation*}
Then 
\begin{equation}
F=\sum_{\{j_{1},\cdots ,j_{k}\}\in \pi _{k}}U^{\alpha _{j_{1}}+\cdots
+\alpha _{j_{k}}}\left( f_{j_{1}}(U^{\alpha _{j_{2}}+\cdots +\alpha
_{j_{k}}}f_{j_{2}}\cdots (U^{\alpha _{j_{k}}}f_{j_{k}})\cdots \right)
\label{re-02}
\end{equation}
where $\pi _{k}$ is the permutation group of $\{1,\cdots ,k\}$.
\end{lemma}

\begin{proof}
We have 
\begin{eqnarray*}
F(x) &=&E^{x}\left\{ \int_{0}^{\infty }\cdots \int_{0}^{\infty }e^{-\alpha
_{1}s_{1}}\cdots e^{-\alpha _{k}s_{k}}f_{1}(X_{s_{1}})\cdots
f_{k}(X_{s_{k}})ds_{1}\cdots ds_{k}\right\} \\
&=&\int_{0}^{\infty }\cdots \int_{0}^{\infty }e^{-\alpha _{1}s_{1}}\cdots
e^{-\alpha _{k}s_{k}}E^{x}\left\{ f_{1}(X_{s_{1}})\cdots
f_{k}(X_{s_{k}})\right\} ds_{1}\cdots ds_{k} \\
&=&\sum_{\{j_{1},\cdots ,j_{k}\}\in \pi _{k}}\idotsint_{0<s_{j_{1}}<\cdots
<s_{j_{k}}<\infty }e^{-\alpha _{1}s_{1}}\cdots e^{-\alpha _{k}s_{k}} \\
&&\times E^{x}\left\{ f_{1}(X_{s_{1}})\cdots f_{k}(X_{s_{k}})\right\}
ds_{1}\cdots ds_{k} \\
&=&\sum_{\{j_{1},\cdots ,j_{k}\}\in \pi _{k}}\idotsint_{0<s_{1}<\cdots
<s_{k}<\infty }e^{-\alpha _{j_{1}}s_{1}}\cdots e^{-\alpha _{j_{k}}s_{k}} \\
&&\times E^{x}\left\{ f_{j_{1}}(X_{s_{1}})\cdots
f_{j_{k}}(X_{s_{k}})\right\} ds_{1}\cdots ds_{k} \\
&=&\sum_{\{j_{1},\cdots ,j_{k}\}\in \pi _{k}}\idotsint_{0<s_{1}<\cdots
<s_{k}<\infty }\int_{E^{\otimes k}}e^{-\alpha _{j_{1}}s_{1}}\cdots
e^{-\alpha _{j_{k}}s_{k}}f_{j_{1}}(z_{1})\cdots f_{j_{k}}(z_{k}) \\
&&\times P_{s_{1}}(x,dz_{1})P_{s_{2}-s_{1}}(z_{1},dz_{2})\cdots
P_{s_{k}-s_{k-1}}(z_{k-1},dz_{k})ds_{1}\cdots ds_{k}
\end{eqnarray*}
and (\ref{re-02}) follows from Lemma \ref{lem1}.
\end{proof}

From now on, we assume that 
\begin{equation*}
X=(\Omega ,\mathscr{F}^{0},\mathscr{F}_{t}^{0},X_{t},\theta _{t},P^{x})
\end{equation*}
is a \emph{continuous Hunt process} in $E^{\prime }=E\cup \{\partial \}$
with the transition semigroup $(P_{t})_{t\geq 0}$. In other words, it is a
Hunt process with continuous sample paths. Therefore, $(X_{t})_{t\geq 0}$ is
a diffusion process in $E$, i.e. $(X_{t})_{t\geq 0}$ possesses the strong
Markov property with continuous sample function. Under our assumptions, any
finite $(\mathscr{F}_{t}^{\mu })$-stopping time is accessible and thus
predictable, and therefore $\mathscr{F}_{T}^{\mu }=\mathscr{F}_{T-}^{\mu }$.
\ In particular, $(\mathscr{F}_{t}^{\mu })$ is left continuous, and thus the
filtration $(\mathscr{F}_{t}^{\mu })$ is continuous for any initial
distribution $\mu $.

Since any martingale on $(\Omega ,\mathscr{F}^{\mu },\mathscr{F}_{t}^{\mu
},P^{\mu })$ has a right continuous modification, by a martingale we always
mean a martingale with right continuous sample function.

\begin{lemma}
\label{l:2.7} Suppose $\xi =\xi _{1}\cdots \xi _{n}$ where each $\xi _{j}$
has the following form 
\begin{equation*}
\xi _{j}=\int_{0}^{\infty }e^{-\alpha _{j}s}f_{j}(X_{s})ds
\end{equation*}
where $\alpha _{j}>0$ and $f_{j}\in \mathscr{B}_{b}(E)$. Let $M_{t}=E^{\mu
}\left\{ \xi |\mathscr{F}_{t}^{\mu }\right\} $. Then $(M_{t})_{t\geq 0}$ is
a bounded continuous martingale on $(\Omega ,\mathscr{F}^{\mu },\mathscr{F}%
_{t}^{\mu },P^{\mu })$.
\end{lemma}

\begin{proof}
According to Lemma~\ref{l:BM}, we need only to show that for function of the
following type 
\begin{equation*}
F(x)= E^{x}\left\{ \left( \prod_{l=k+1}^{n}e^{-\alpha
_{j_{l}}t}\int_{0}^{\infty }e^{-\alpha _{j_{l}}s}f_{j_{l}}(X_{s})ds\right)
\right\},
\end{equation*}
$t\mapsto F(X_t)$ is continuous. By Lemma~\ref{lem2}, $F$ is an $\alpha $%
-potential, so that it is finely continuous, and together with Lemma \ref
{meyer-lemma}, it implies that $t\rightarrow F(X_t)$ is continuous, which
completes the proof.
\end{proof}

We now state the main result of this section. For simplicity, a square
integrable martingale $(M_t)_{t\ge 0}$ over $(\Omega,\mathscr{M},\mathscr{M}%
_t,P)$ means $M_t=E(\xi|\mathscr{M}_t)$ with $\xi\in L^2(\Omega,\mathscr{M}%
,P)$. This is equivalent to say $\sup_{t>0}E[M^2_t]<\infty$.

\begin{theorem}
\label{mc-th}Let 
\begin{equation*}
X=(\Omega ,\mathscr{F}^{0},\mathscr{F}_{t}^{0},X_{t},\theta _{t},P^{x})
\end{equation*}
be a continuous Hunt process in $E$, and $\mu \in \mathscr{P}(E)$. If $\xi
\in L^{2}(\Omega ,\mathscr{F}^{\mu },P^{\mu })$, then the martingale $%
M_{t}=E^{\mu }\left\{ \xi |\mathscr{F}_{t}^{\mu }\right\} $ is continuous,
that is, square-integrable martingales on $(\Omega ,\mathscr{F}^{\mu },%
\mathscr{F}_{t}^{\mu },P^{\mu })$ are continuous. Therefore local
martingales on $(\Omega ,\mathscr{F}^{\mu },\mathscr{F}_{t}^{\mu },P^{\mu })$
are continuous.
\end{theorem}

\begin{proof}
We can choose a sequence $\xi _{n}\in \mathscr{C}$ such that $\xi
_{n}\rightarrow \xi $ in $L^{2}$. Doob's maximal inequality implies that, if
necessary by considering a subsequence, the martingales $\{E^{\mu }(\xi _{n}|%
\mathscr{F}_{t}^{\mu }):t\ge 0\}$ converges (almost surely at least along a
subsequence) to $\{E^{\mu }(\xi |\mathscr{F}_{t}^{\mu }):t\ge 0\}$ uniformly
on any finite interval of $t\ge 0$. It is shown in Lemma~\ref{l:2.7} that
for each $n$, the martingale $E^{\mu }(\xi _{n}|\mathscr{F}_{t}^{\mu })$ is
continuous and thus the square integrable martingale $\{E^{\mu }(\xi |%
\mathscr{F}_{t}^{\mu }):t\ge 0\}$ must be continuous.

By the localization technique, it follows thus that local martingales on $%
(\Omega ,\mathscr{F}^{\mu },\mathscr{F}_{t}^{\mu },P^{\mu })$ are continuous.
\end{proof}

\section{Martingale representation for continuous Hunt process}

In this section we assume that 
\begin{equation*}
X=(\Omega ,\mathscr{F}^{0},\mathscr{F}_{t}^{0},X_{t},\theta _{t},P^{x})
\end{equation*}
is a \emph{continuous} Hunt process in the state space $E^{\prime }=E\cup
\{\partial \}$ with transition semigroup $\{P_{t}(x,dy):t\geq 0\}$, where $E$
is a locally compact separable metric space. Let $\mu \in \mathscr{P}(E)$ be
an initial distribution.\ 

If $\alpha >0$ and $f\in \mathscr{B}_{b}(E)$ then $M^{\alpha ,f}$ denotes
the continuous martingale 
\begin{equation*}
M_{t}^{\alpha ,f}=E^{\mu }\left\{ \int_{0}^{\infty }e^{-\alpha s}f(X_{s})ds|%
\mathscr{F}_{t}^{\mu }\right\} \text{.}
\end{equation*}
Recall that, if $u$ is an $\alpha $-potential, i.e., $u=U^{\alpha }f$ where $%
f\in \mathscr{B}_{b}(E)$, then $u(X_{t})-u(X_{0})$ is a continuous
semimartingale on $(\mathscr{F}^{\mu },\mathscr{F}_{t}^{\mu },P^{\mu })$. 
% with continuous sample function (more precisely, its
%sample function has continuous modification -- this comment applies to
%similar situation if no confusion may arise)
and possesses Doob-Meyer's decomposition 
\begin{equation*}
u(X_{t})-u(X_{0})=M_{t}^{[u]}+A_{t}^{[u]}
\end{equation*}
where 
\begin{equation*}
M_{t}^{[u]}=\int_{0}^{t}e^{\alpha s}dM_{s}^{\alpha ,f}\text{, \ }%
A_{t}^{[u]}=\int_{0}^{t}Lu(X_{s})ds
\end{equation*}
and $Lu=\alpha u-f$. 
%For simplicity, such martingale $M_{t}^{[u]}$ is called a resolvent martingale.

We make the following assumptions on the continuous Hunt process $X$ started
with an initial distribution $\mu \in \mathscr{P}(E)$, and we call these
assumptions the\emph{\ Fukushima representation property}.

\noindent \textbf{Assumptions}. There is an algebra (a vector space which is
closed under the multiplication of functions) $K(E)\subset \mathscr{B}%
_{b}(E) $ which generates the Borel $\sigma $-algebra $\mathscr{B}(E)$ and
is invariant under $U^{\alpha}$ for $\alpha>0$, and there are finite many
continuous martingales $M^{1},\cdots ,M^{d}$ over $(\Omega ,\mathscr{F}^{\mu
},\mathscr{F}_{t}^{\mu },P^{\mu })$ such that the following conditions are
satisfied:

\begin{itemize}
\item[(1)]  For any potential $u=U^{\alpha }f$ where $\alpha >0$ and $f\in
K(E)$, the martingale part $M^{[u]}$ of the semimartingale $%
u(X_{t})-u(X_{0}) $ has the martingale representation in terms of $%
(M^{1},\cdots ,M^{d})$, that is, there are predictable processes $F_{1}$, $%
\cdots $, $F_{d}$ on $(\Omega ,\mathscr{F}^{\mu },\mathscr{F}_{t}^{\mu })$
such that 
\begin{equation}
M_{t}^{[u]}=\sum_{j=1}^{d}\int_{0}^{t}F_{s}^{j}dM_{s}^{j}\text{ \ \ }P^{\mu }%
\text{-a.e.}  \label{kr1}
\end{equation}

\item[(2)]  $\left( \langle M^{j},M^{i}\rangle _{t}\right) $ is strictly
positive definite.
\end{itemize}

The first assumption means that the martingale $M^{[u]}$ with $u$ being a
potential may be represented. The second condition ensures that the
representation (\ref{kr1}) is unique.

The Fukushima representation property is mainly an abstraction of the chain
role for the martingale part of $u(X_{t})$. Indeed, if $X_{t}=(X_{t}^{1},%
\cdots ,X_{t}^{d})$ is a $d$-dimensional Brownian motion and $u$ is an $%
\alpha$-potential with $\alpha>0$, then $u$ is smooth and by It\^{o}'s
formula 
\begin{equation*}
u(X_{t})-u(X_{0})=\sum_{j=1}^{d}\int_{0}^{t}\frac{\partial u}{\partial x^{j}}%
(X_{s})dX_{s}^{j}+\int_{0}^{t}\frac{1}{2}\Delta u(X_{s})ds
\end{equation*}
so that 
\begin{equation*}
M_{t}^{[u]}=\sum_{j=1}^{d}\int_{0}^{t}\frac{\partial u}{\partial x^{j}}%
(X_{s})dX_{s}^{j}\text{.}
\end{equation*}
One can easily see that the Brownian motion satisfies the Fukushima
representation property.

\begin{theorem}[Martingale representation]
\label{th1-d} Let $\mu \in \mathscr{P}(E)$. Suppose that the Fukushima
representation property holds for $X$ with a finite set of martingales $%
(M^{1},\cdots ,M^{d})$. For any square-integrable martingale $%
N=(N_{t})_{t\geq 0}$ on $(\Omega ,\mathscr{F}^{\mu },\mathscr{F}_{t}^{\mu
},P^{\mu })$, there are unique predictable processes $(F_{t}^{i})$ such that 
\begin{equation*}
N_{t}-N_{0}=\sum_{i=1}^{d}\int_{0}^{t}F_{s}^{i}dM_{s}^{i}\text{ \ \ \ \ }%
P^{\mu }\text{-a.e.}
\end{equation*}
\end{theorem}

\begin{proof}
The uniqueness follows from condition 2) in the Fukushima representation. We
prove the existence. Take $\xi \in L^{2}(\Omega ,\mathscr{F}^{\mu },P^{\mu
}) $ such that $N_{t}=E^{\mu }\{\xi |\mathscr{F}_{t}^{\mu }\}$. Since $%
\mathscr{C}$ is dense in $L^{2}(\Omega ,\mathscr{F}^{\mu },P^{\mu })$, so we
first prove the martingale representation for $\xi \in $ $\mathscr{C}$. By
the linearity, we only need to consider the case that $\xi =\xi _{1}\cdots
\xi _{n}$ where $\xi _{j}=\int_{0}^{\infty }e^{-\alpha _{j}s}f_{j}(X_{s})ds$
for $\alpha _{j}>0$ and $f_{j}\in K(E)$. In this case, according to \ref
{jk-01}, Lemma \ref{lem1} and Lemma \ref{lem2} 
\begin{equation*}
N_{t}=E^{\mu }\{\xi |\mathscr{F}_{t}^{\mu }\}=\sum_{m}Z_{t}^{m}
\end{equation*}
where the sum is a finite one, and for each $m$, $Z^{m}=Z_{t}$ has the
following form 
\begin{equation*}
Z_{t}^{m}=V_{t}^{m}u^{m}(X_{t})
\end{equation*}
(the superscript $m$ will be dropped if no confusion may arise), where 
\begin{equation*}
V_{t}=\prod_{i=1}^{k^{\prime }}\int_{0}^{t}e^{-\beta _{i}s}g_{i}(X_{s})ds
\end{equation*}
and 
\begin{eqnarray*}
u(x) &=&\idotsint_{0<s_{1}<\cdots <s_{k}<\infty }e^{-\sum_{j=1}^{k}\beta
_{j}s_{j}}\int_{E^{\otimes k}}h_{1}(z_{1})\cdots
h_{k}(z_{k})P_{s_{1}}(x,dz_{1}) \\
&&\times P_{s_{2}-s_{1}}(z_{1},dz_{2})\cdots
P_{s_{k}-s_{k-1}}(z_{k-1},dz_{k})ds_{1}\cdots ds_{k}
\end{eqnarray*}
for some $k^{\prime }$ and $k$, $\beta _{i}>0$ and functions $g_{i}$, $h_{j}$
are bounded and continuous. According to Lemma~\ref{lem1} 
\begin{equation*}
u=U^{\beta _{1}+\cdots +\beta _{k}}\left( h_{1}(U^{\beta _{2}+\cdots +\beta
_{k}}h_{2}\cdots (U^{\beta _{k}}h_{k})\cdots \right) \text{.}
\end{equation*}
In particular, $u$ is again a potential which has a form $u=U^{\alpha }g$
for 
\begin{equation*}
g=h_{1}(U^{\beta _{2}+\cdots +\beta _{k}}h_{2}\cdots (U^{\beta
_{k}}h_{k})\cdots )\in K(E)
\end{equation*}
and $\alpha =\beta _{1}+\cdots +\beta _{k}$. Hence $u(X_{t})$ is a
continuous semimartingale with decomposition 
\begin{equation*}
u(X_{t})-u(X_{0})=M_{t}^{[u]}+A_{t}^{[u]}
\end{equation*}
where $A^{[u]}$ is continuous with finite variation, and due to the
Fukushima representation property 
\begin{equation*}
M_{t}^{[u]}=\sum_{j=1}^{d}\int_{0}^{t}G_{s}^{j}dM_{s}^{j}
\end{equation*}
for some predictable processes $G^{j}$. In particular, each $Z^{m}$ is a
continuous semimartingale. Since, by Theorem \ref{mc-th}, $N$ is a
continuous martingale, so that 
\begin{equation*}
N_{t}=\sum_{m}\text{ the continuous martingale part of }V_{t}^{m}u^{m}(X_{t})%
\text{.}
\end{equation*}
Therefore we are interested in the martingale part of $Z_{t}=V_{t}u(X_{t})$.
Since $V$ is a\ finite variation process, so according to It\^{o}'s formula 
\begin{eqnarray*}
Z_{t} &=&Z_{0}+\int_{0}^{t}u(X_{s})dV_{s}+\int_{0}^{t}V_{s}du(X_{s}) \\
&=&Z_{0}+\int_{0}^{t}u(X_{s})dV_{s}+\int_{0}^{t}V_{s}dA_{t}^{[u]}+%
\int_{0}^{t}V_{s}dM_{t}^{[u]} \\
&=&Z_{0}+\int_{0}^{t}u(X_{s})dV_{s}+\int_{0}^{t}V_{s}dA_{t}^{[u]}+%
\sum_{i=1}^{d}\int_{0}^{t}V_{s}\cdot G_{s}^{i}dM_{s}^{i}
\end{eqnarray*}
so that the martingale part of $Z_{t}$ is 
\begin{equation*}
\sum_{i=1}^{d}\int_{0}^{t}V_{s}\cdot G_{s}^{i}dM_{s}^{i}\text{.}
\end{equation*}
Therefore 
\begin{equation*}
N_{t}=E^{\mu }\{\xi |\mathscr{F}_{t}^{\mu
}\}=\sum_{i=1}^{d}\int_{0}^{t}\sum_{m}V_{s}^{m}\cdot G_{s}^{m,i}dM_{s}^{i}
\end{equation*}
which shows the martingale representation.

Suppose now $\xi \in L^{2}(\Omega ,\mathscr{F}^{\mu },P^{\mu })$. Choose a
sequence $\xi _{n}\in \mathscr{C}$ such that $\xi _{n}\rightarrow \xi $ in $%
L^{2}(\Omega ,\mathscr{F}^{\mu },P^{\mu })$. Let $N_{t}^{(n)}=E^{\mu }(\xi
_{n}|\mathscr{F}_{t}^{\mu })$ and $N_{t}=E^{\mu }(\xi |\mathscr{F}_{t}^{\mu
})$. According to Doob's maximal inequality, if necessary by passing to a
subsequence, we can assume that $N_{t}^{(n)}$ converges to $N_{t}$ uniformly
on any finite interval. $N_{t}^{(n)}$ has the martingale representation 
\begin{equation*}
N_{t}^{(n)}-N_{0}^{(n)}=\sum_{j=1}^{d}\int_{0}^{t}F(n)_{s}^{j}dM_{s}^{j}
\end{equation*}
so that 
\begin{eqnarray*}
&&\langle N_{t}^{(n)}-N_{t}^{(m)},N_{t}^{(n)}-N_{t}^{(m)}\rangle \\
&=&\sum_{i,j=1}^{d}%
\int_{0}^{t}(F(n)_{s}^{i}-F(m)_{s}^{i})(F(n)_{s}^{j}-F(m)_{s}^{j})d\langle
M^{i},M^{j}\rangle _{s}\text{.}
\end{eqnarray*}
Since $(\langle M^{i},M^{j}\rangle _{t})$ is positive, it follows that $%
(F(n)^{1},\cdots ,F(n)^{d})$ converges to predictable processes $%
(F^{1},\cdots ,F^{d})$ under the norm 
\begin{equation}
||(F^{1},\cdots ,F^{d})||=\sum_{N=1}^{\infty }\frac{1}{2^{N}}%
\sum_{i,j=1}^{d}E^{\mu }\left[ \int_{0}^{N}F_{s}^{i}F_{s}^{j}d\langle
M^{i},M^{j}\rangle _{s}\right] \text{.}  \label{norm-a}
\end{equation}
Then 
\begin{equation*}
N_{t}-N_{0}=\sum_{j=1}^{d}\int_{0}^{t}F_{s}^{j}dM_{s}^{j}\text{.}
\end{equation*}
\end{proof}

This theorem claims that as long as every martingale of resolvent type is
representable, so is any martingale. When is the Fukushima representation
property satisfied? There are many examples. In the remain of this section,
we shall give three interesting examples in symmetric situation.

Brownian motion with any initial distribution is certainly an example.
Indeed, for Brownian motion in $\mathbb{R}^{d}$, we may choose $%
K(E)=C_{\infty }^{\infty }(R^{d})$ (the space of smooth functions which
vanish at infinity), then for $f\in K(E)$, $U^{\alpha }f$ is smooth, and (%
\ref{kr1}) follows from It\^{o}'s formula applying to $U^{\alpha }f$.
Theorem~\ref{th1-d} gives a new proof for classical martingale
representation theorem.

The second example is the reflecting Brownian motion. As Example 1.6.1 in 
\cite{FOT}, we consider Dirichlet form $({\frac{1}{2}}\mathbf{D},H^{1}(D))$
on $L^{2}(D)$ where $\mathbf{D}$ is the classical Dirichlet integral and $D$
is a bounded domain on $\mathbb{R}^{d}$. We further assume that any $x\in
\partial D$ has a neighborhood $U$ such that 
\begin{equation*}
D\cap U=\{(x_{i})\in \mathbb{R}^{d}:x_{d}>F(x_{1},\cdots ,x_{d-1})\}\cap U
\end{equation*}
for some continuous function $F$. Then $C_{0}^{\infty }(\overline{D})$ (the
space of restriction to $\overline{D}$ of functions in $C_{0}^{\infty }(%
\mathbb{R}^{d})$) is dense in $H^{1}(D)$ (see \cite{mazja} for details),
i.e., $({\frac{1}{2}}\mathbf{D},H^{1}(D))$ is a regular Dirichlet form on $%
L^{2}(\overline{D})$. The corresponding continuous Hunt process $%
X=(X_{t},P^{x})$ is called the reflecting Brownian motion. For $x=(x^{i})\in 
\mathbb{R}^{d}$, we use $u_{i}(x)=x^{i}$, $1\leq i\leq d$, to denote the
coordinate functions. Then $u_{i}\in \mathscr{F}$ and we denote by $%
M^{i}=M^{[u_{i}]}$ the martingale part in Fukushima's decomposition. It can
be seen from Corollary 5.6.2 \cite{FOT} that for any $u\in C_{0}^{\infty }(%
\overline{D})$, 
\begin{equation*}
M_{t}^{[u]}=\sum_{i=1}^{d}\int_{0}^{t}{\frac{\partial u}{\partial x_{i}}}%
(X_{s})dM_{s}^{i},\ P^{x}\text{-a.s. for q.e.}\ x\in \overline{D},
\end{equation*}
where q.e. means `quasi-everywhere', i.e., except a set of zero-capacity.
Then a routine approximation procedure shows that for any $u\in H^{1}(D)$,
there exist Borel measurable functions $\{f_{i}:1\leq i\leq d\}$ on $%
\overline{D}$ such that 
\begin{equation*}
M_{t}^{[u]}=\sum_{i=1}^{d}\int_{0}^{t}f_{i}(X_{s})dM_{s}^{i}\text{,}\ P^{x}%
\text{-a.s. for q.e.}\ x\in \mathbb{R}^{d}.
\end{equation*}
Therefore the reflecting Brownian motion has Fukushima representation
property, by choosing $K(\overline{D})$ to be the space of bounded
measurable functions and any initial distribution $\mu $ charging no set of
zero capacity, i.e., a smooth distribution, because an exceptional set
exists in above representation as is always when the process is constructed
through a Dirichlet form. If the boundary is Lipschitz, then the transition
function has density (\cite{bass1}) and in this case, the exceptional set
may be erased. Notice that under the current condition, the reflecting
Brownian motion $X$ itself is not necessarily a semimartingale. The readers
who are interested may refer to \cite{bass1}, \cite{Chen1} and \cite{Chen2}
about when a reflecting BM is a semimartingale and the corresponding
Skorohod decomposition. It should be pointed out that, although the
martingale part of the reflected Brownian motion is a Brownian motion, but
the martingale representation property does not follow from the classical
representation property for Brownian motion. The reason is that, as long as
the boundary is not sufficiently smooth, the natural filtration $(\mathscr{F}%
_{t}^{\mu })_{t\geq 0}$ is much bigger in general than the natural
filtration generated by the martingale part $(M^{1},\cdots ,M^{d})$ of $X$.

Another example our main result may apply is symmetric diffusions in a
domain killed at boundary. Actually Theorem 6.2.2 in \cite{FOT} tells us
that every continuous symmetric Hunt process with a smooth core enjoys\ the
Fukushima representation property. More precisely let $D$ be a domain of $%
\mathbb{R}^{d}$ with continuous boundary $\partial D$ and $m$ a Radon
measure on $D$. Let $X$ be a continuous Hunt process which is symmetric with
respect to $m$ and $(\mathscr{E},\mathscr{F})$ the associated Dirichlet form
on $L^{2}(D,m)$, which has $C_{0}^{1}(D)$ as a core. For $x=(x^{i})\in 
\mathbb{R}^{d}$, we use $u_{i}(x)=x^{i}$, $1\leq i\leq d$, to denote the
coordinate functions. Then $u_{i}\in \mathscr{F}_{\mathrm{loc}}$ and we
denote by $M^{i}=M^{[u_{i}]}$ the martingale part in Fukushima's
decomposition. Let 
\begin{equation*}
\mu _{i,j}=\mu _{\ang{M^i,M^j}},\ 1\leq i,j\leq d,
\end{equation*}
the smooth measure associated with CAF $\ang{M^i,M^j}$. Then $\mathscr{E}$
is expressed as 
\begin{equation*}
\mathscr{E}(u,v)=\sum_{i,j=1}^{d}\int_{D}{\frac{\partial u}{\partial x^{i}}}{%
\frac{\partial u}{\partial x^{j}}}d\mu _{i,j}(x),\ u,v\in C_{0}^{1}(D)\text{.%
}
\end{equation*}

As asserted in Theorem 6.2.2 \cite{FOT}, for any initial smooth distribution 
$\mu $ (i.e. a probability on $(D,\mathscr{B}(D))$ having no charge on
capacity zero sets) and $u\in \mathscr{F}$, the martingale part $M^{[u]}$ in
Fukushima's decomposition of $u$ may be represented as 
\begin{equation*}
M_{t}^{[u]}=\sum_{i=1}^{d}\int_{0}^{t}f_{i}(X_{s})dM_{s}^{i}\text{ \ \ \ \ \
\ }P^{\mu }\text{-a.e.}
\end{equation*}
where $f_{1},\cdots ,f_{d}\in \mathscr{B}(D)$. If we take $%
K(E)=L^{2}(E,m)\cap \mathscr{B}_{b}(D)$, $X$ satisfies the Fukushima
representation property. In these examples, $\{M^{i}\}$ are the martingales
corresponding to coordinate functions so we call them coordinate martingales.

To have the uniqueness, some kind of non-degenerateness is needed. We say
that $X$ is non-degenerate if the condition (2) in Fukushima representation
property is satisfied: $(\ang{M^i,M^j})_{1\le i,j\le d}$ is positive.

\begin{corollary}
\label{coro1} Assume that $X$ is either the reflecting Brownian motion on a
bounded domain or a non-degenerate symmetric Hunt diffusion on a domain $%
D\subset \mathbb{R}^{d}$ as stated above. Then the Fukushima representation
property is satisfied and therefore the martingale representation holds in
the sense of Theorem~\ref{th1-d} with coordinate martingales and for a given
initial distribution $\mu$ charging no sets of zero capacity.
\end{corollary}

From this result, we may recover the martingale representation established
in \cite{bally} and \cite{zheng2}, where $X$ is a diffusion process
corresponding to non-degenerate symmetric elliptic operator on $\mathbb{R}%
^{d}$.

Without essential difference, the conclusion holds also for reflecting
diffusions on such domain with generator being a symmetric uniformly
elliptic differential operator of second order as introduced in the
beginning of next section.

\section{Backward stochastic differential equations}

In this section we consider backward stochastic differential equations which
can be used to provide probability representations for weak solutions of the
initial and boundary value problem of a quasi-linear parabolic equation.

Let $D\subset \mathbb{R}^{d}$ be a bounded domain with a continuous boundary 
$\partial D$, $\overline{D}=D\cup \partial D$ the closure of $D$. Let 
\begin{equation*}
L=\frac{1}{2}\sum_{i,j=1}^{d}\frac{\partial }{\partial x^{j}}a^{ij}\frac{%
\partial }{\partial x^{i}}\text{ }
\end{equation*}
be an elliptic differential operator of second order, where $a=(a^{ij})$ is
a positive-definite, symmetric, matrix-valued function on $D$, $a=(a^{ij})$
is Borel measurable, and satisfies the elliptic condition: 
\begin{equation*}
\lambda |\xi |^{2}\leq \sum_{i,j=1}^{d}a^{ij}(x)\xi _{i}\xi _{j}\leq \lambda
^{-1}|\xi |^{2}\text{ \ \ \ \ \ }\forall \xi =(\xi _{i})\in \mathbb{R}^{d}
\end{equation*}
for all $x\in D$ for some constant $\lambda >0$. Consider the Dirichlet form 
$(\mathscr{E},\mathscr{F})$ on $L^{2}(D,dx)$, where 
\begin{equation}
\mathscr{E}(u,v)=\frac{1}{2}\int_{D}\sum_{i,j=1}^{d}a^{ij}\frac{\partial u}{%
\partial x^{j}}\frac{\partial v}{\partial x^{i}}  \label{dr-1}
\end{equation}
and $\mathscr{F}=H^{1}(D)$.

Let $\Omega $ be a space of all continuous paths in $\overline{D}$, $%
(X_{t})_{t\geq 0}$ the coordinate process on $\Omega $, $\mathscr{F}%
^{0}=\sigma \{X_{s}:s\geq 0\}$, $\mathscr{F}_{t}^{0}=\sigma \{X_{s}:s\leq
t\} $ for each $t\geq 0$, and $(\theta _{t})_{t\geq 0}$ shift operators on $%
\Omega $. Let 
\begin{equation*}
X=(\Omega ,\mathscr{F}^{0},\mathscr{F}_{t}^{0},X_{t},\theta _{t},P^{x})
\end{equation*}
be the canonical realization of the symmetric diffusion process in the state
space $\overline{D}$ associated with the Dirichlet space $(\mathscr{E},%
\mathscr{F})$, which is called a reflecting symmetric diffusion in $D$.

The coordinate functions $u_{j}(x)=x^{j}$ ($j=1,\cdots ,d$) belong to the
local Dirichlet space $\mathscr{F}_{\text{loc}}$, so that 
\begin{equation}  \label{e:fuku1}
X_{t}^{j}-X_{0}^{j}=M_{t}^{j}+A_{t}^{j}\text{ \ \ \ }P^{x}\text{-a.e.\ \ }%
j=1,\cdots ,d
\end{equation}
for all $x\in \overline{D}$ except for a capacity zero set, where $%
M^{j}=M^{[f_{j}]}$ etc.

Let $\mathscr{S}_{1}(\overline{D})$ denote the space of all probability $\mu
\in \mathscr{P}(\overline{D})$ which has no charge on zero capacity sets
(with respect to the Dirichlet form $(\mathscr{E},H^{1}(D))$ defined by (\ref
{dr-1}). According to Theorem \ref{th1-d}, for any initial distribution $\mu
\in \mathscr{S}_{1}(\overline{D})$, \ the family of martingales\ $%
\{M^{j}:j=1,\cdots ,d\}$ over $(\Omega ,\mathscr{F}^{\mu },\mathscr{F}%
_{t}^{\mu },P^{\mu })$ has the martingale representation property: for any
square-integrable martingale $N=(N_{t})_{t\geq 0}$ on $(\Omega ,\mathscr{F}%
^{\mu },\mathscr{F}_{t}^{\mu },P^{\mu })$, there are unique predictable
processes $(F_{t}^{i})$ such that 
\begin{equation*}
N_{t}-N_{0}=\sum_{i=1}^{d}\int_{0}^{t}F_{s}^{i}dM_{s}^{i}\text{ \ \ \ \ }%
P^{\mu }\text{-a.e.}
\end{equation*}

Let us work with a fixed smooth initial distribution $\mu \in \mathscr{S}%
_{1}(D)$ and the filtered probability space $(\Omega ,\mathscr{F}^{\mu },%
\mathscr{F}_{t}^{\mu },P^{\mu })$.

Consider the following backward stochastic differential equation 
\begin{equation}
dY_{t}^{i}=-f^{i}(t,Y_{t},Z_{t})dt+\sum_{i,j=1}^{d}Z_{t}^{ij}dM_{t}^{j}\text{%
, \ }Y_{T}^{i}=\xi ^{i}  \label{bsde=t1}
\end{equation}
$i=1,\cdots ,d^{\prime }$, where $T>0$, $\xi ^{i}\in L^{2}(\Omega ,%
\mathscr{F}_{T}^{\mu },P^{\mu })$ are given terminal values, and $f^{i}$ are
Lipschitz functions: there is a constant $C_{1}\geq 0$%
\begin{equation*}
|f^{i}(t,y,z)|\leq C_{1}(1+t+|y|+|z|)
\end{equation*}
and 
\begin{equation*}
|f^{i}(t,y,z)-f^{i}(t,\tilde{y},\tilde{z})|\leq C_{1}\left( |y-\tilde{y}|+|z-%
\tilde{z}|\right) 
\end{equation*}
for all $t\geq 0$, $y,\tilde{y}\in \mathbb{R}^{d^{\prime }}$, $z,\tilde{z}%
\in \mathbb{R}^{d^{\prime }\times d}$. One seeks for a solution pair $(Y,Z)$
which solves the following integral equation 
\begin{equation}
Y_{t}^{i}-\xi
^{i}=\int_{t}^{T}f^{i}(s,Y_{s},Z_{s})ds-\sum_{j=1}^{d}%
\int_{t}^{T}Z_{s}^{ij}dM_{s}^{j}\text{ }  \label{bsde=t2}
\end{equation}
for $t\in \lbrack 0,T]$. The integral equation (\ref{bsde=t2}) has a unique
solution pair $(Y,Z)$ such that $Y^{i}$ is a continuous semimartingale, and $%
Z^{ij}$ are predictable processes satisfying 
\begin{equation*}
E^{\mu
}\int_{0}^{T}\sum_{k,l=1}^{d}a^{kl}(X_{s})Z_{s}^{il}Z_{s}^{ki}ds<\infty 
\text{.}
\end{equation*}
This can be demonstrated by employing the Picard iteration for $(Y,Z)$ as in
the case of Brownian motion (see \cite{pardoux1}). Another approach,
proposed in a paper by Lyons, Liang and Qian \cite{lyons-etc1} which applies
to a general filtered probability space, may be described as follows. The
idea is to rewrite the integral equation (\ref{bsde=t2}) into a functional
differential equation for the variation process part $V$ of $Y$. Let $Y=N-V$
where $V$ is a finite variation process, and 
\begin{equation*}
N_{t}^{i}-N_{0}^{i}=\sum_{j=1}^{d}\int_{0}^{t}Z_{s}^{ij}dM_{s}^{j}\text{.}
\end{equation*}
On the other hand 
\begin{equation*}
N_{t}=E^{\mu }\{\xi +V_{T}|\mathscr{F}_{t}^{\mu }\}\text{.}
\end{equation*}
Since $Y$ is a continuous semimartingale, its decomposition is unique up to
an initial value. The integral equation (\ref{bsde=t2}) leads to that 
\begin{equation*}
V_{t}=-\int_{t}^{T}f^{i}(s,Y_{s},Z_{s})ds+N_{T}-\xi 
\end{equation*}
conditioned on $\mathscr{F}_{t}^{\mu }$ and we obtain 
\begin{eqnarray*}
V_{t} &=&-E^{\mu }\left\{ \int_{t}^{T}f^{i}(s,Y_{s},Z_{s})ds|\mathscr{F}%
_{t}^{\mu }\right\} +N_{t}-E^{\mu }\left\{ \xi |\mathscr{F}_{t}^{\mu
}\right\}  \\
&=&-E^{\mu }\left\{ \int_{t}^{T}f^{i}(s,Y_{s},Z_{s})ds|\mathscr{F}_{t}^{\mu
}\right\} +E^{\mu }\left\{ V_{T}|\mathscr{F}_{t}^{\mu }\right\}  \\
&=&-E^{\mu }\left\{ \int_{0}^{T}f^{i}(s,Y_{s},Z_{s})ds-V_{T}|\mathscr{F}%
_{t}^{\mu }\right\} +\int_{0}^{t}f^{i}(s,Y_{s},Z_{s})ds.
\end{eqnarray*}
Therefore the integral equation (\ref{bsde=t2}) is equivalent to 
\begin{equation}
V_{t}-V_{0}=\int_{0}^{t}f^{i}(s,Y_{s},Z_{s})ds  \label{zs-1}
\end{equation}
where 
\begin{equation*}
Y_{t}=Y(V)_{t}=N(V)_{t}-V_{t}\text{, \ \ }N(V)_{t}=E^{\mu }\left\{ \xi
+V_{T}|\mathscr{F}_{t}^{\mu }\right\} 
\end{equation*}
and $Z_{t}=Z(V)_{t}$ is determined by the martingale representation theorem 
\begin{equation*}
N(V)_{t}^{i}-N(V)_{0}^{i}=\sum_{j=1}^{d}\int_{0}^{t}Z(V)_{s}^{ij}dM_{s}^{j}%
\text{.}
\end{equation*}
Equation (\ref{zs-1}) thus may be written as a functional equation 
\begin{equation}
V_{t}-V_{0}=\int_{0}^{t}f^{i}(s,Y(V)_{s},Z(V)_{s})ds  \label{zs-2}
\end{equation}
where $Y(V)$ and $Z(V)$ are considered as functionals of $V$. The Picard
iteration applies to (\ref{zs-2}) we have

\begin{theorem}
\label{bsder}If $\xi \in L^{2}(\Omega ,\mathscr{F}_{T}^{\mu },P^{\mu })$ and 
$f^{i}$ are Lipschitz continuous, then there is a unique pair $(Y,Z)$ such
that $Y$ is a continuous semimartingale which solves BSDE (\ref{bsde=t1}).
\end{theorem}

For a complete proof of Theorem \ref{bsder}, the reader may refer to \cite
{lyons-etc1}.

\section{Non-linear parabolic equations}

We are under the same setting as in the previous section, and use the
notations established therein.

To motivate our approach, let us begin with the case that $a$ is smooth, and 
$D$ is bounded domain with a smooth boundary.

In this case $X^{j}=(X_{t}^{j})_{t\geq 0}$ in \eqref{e:fuku1} are continuous
semimartingales, thus $A^{j}$ are finite variation processes. For any $h\in
C_{b}^{1,2}([0,\infty )\times \overline{D})$ satisfying the Neumann boundary
condition that $\left. \frac{\partial h}{\partial \nu }\right| _{\partial
D}=0$, where $\frac{\partial }{\partial \nu }$ denotes the normal derivative
with respect to the Riemann metric $(a^{ij})=(a_{ij})^{-1}$, we have 
\begin{equation*}
h(t,X_{t})-h(0,X_{0})=M_{t}^{h}+A_{t}^{h}
\end{equation*}
where 
\begin{equation}
M_{t}^{h}=h(t,X_{t})-h(0,X_{0})-\int_{0}^{t}\left( \frac{\partial }{\partial
s}+L\right) h(s,X_{s})ds  \label{m-1}
\end{equation}
is a martingale under $P^{x}$, and 
\begin{equation}
A_{t}^{h}=\int_{0}^{t}\left( \frac{\partial }{\partial s}+L\right)
h(s,X_{s})ds\text{.}  \label{m-3-a}
\end{equation}

On the other hand, applying It\^{o}'s formula 
\begin{eqnarray*}
h(t,X_{t})-h(0,X_{0}) &=&\int_{0}^{t}\left( \frac{\partial }{\partial s}+%
\frac{1}{2}\sum_{i,j=1}^{d}a^{ij}\frac{\partial ^{2}}{\partial x^{i}\partial
x^{j}}\right) h(s,X_{s})ds \\
&&+\int_{0}^{t}\sum_{j=1}^{d}\frac{\partial }{\partial x^{j}}%
h(s,X_{s})d(M_{s}^{j}+A_{s}^{j})
\end{eqnarray*}
it thus follows that 
\begin{equation}
M_{t}^{h}=\sum_{j=1}^{d}\int_{0}^{t}\frac{\partial }{\partial x^{j}}%
h(s,X_{s})dM_{s}^{j}  \label{m-2}
\end{equation}
and 
\begin{equation}
A_{t}^{i}=\frac{1}{2}\int_{0}^{t}\sum_{j=1}^{d}\frac{\partial }{\partial
x^{j}}a^{ij}(X_{s})ds\text{.}  \label{m-3}
\end{equation}

Consider a solution $u(x,t)$ to the initial boundary problem to the
non-linear parabolic equation 
\begin{equation}
\begin{cases}
\displaystyle{\left( \frac{\partial }{\partial t}-L\right) u+f(t,u,\nabla
u)=0}\text{,} & \ \label{no-l1} \\ 
u(0,x)=\varphi (x), & x\in \mathbb{R}^{d}, \\ 
\displaystyle{\ \frac{\partial u(t,\cdot )}{\partial \nu }\biggm|_{\partial
D}=0}, & t>0
\end{cases}
\end{equation}
Then, by (\ref{m-1}) and (\ref{m-2}) 
\begin{eqnarray*}
h(T,X_{T}) &=&h(t,X_{t})+\int_{t}^{T}\left( \frac{\partial }{\partial s}%
+L\right) h(s,X_{s})ds \\
&&+\sum_{j=1}^{d}\int_{t}^{T}\frac{\partial }{\partial x^{j}}%
h(s,X_{s})dM_{s}^{j}
\end{eqnarray*}
together with the PDE (\ref{no-l1}) we deduce that 
\begin{eqnarray}
h(t,X_{t})-h(T,X_{T}) &=&\int_{t}^{T}f(T-s,h(s,X_{s}),\nabla h(s,X_{s}))ds 
\notag \\
&&-\sum_{j=1}^{d}\int_{t}^{T}\frac{\partial }{\partial x^{j}}%
h(s,X_{s})dM_{s}^{j}\text{.}  \label{fz1}
\end{eqnarray}

Let $Y_{t}=u(T-t,X_{t})$ and $Z_{t}^{j}=\frac{\partial }{\partial x^{j}}%
h(t,X_{t})$. Then the previous equation may be written as 
\begin{equation}
Y_{t}-Y_{T}=\int_{t}^{T}f(T-s,Y_{s},Z_{s})ds-\sum_{j=1}^{d}%
\int_{t}^{T}Z_{s}^{j}dM_{s}^{j}  \label{bs-1}
\end{equation}
and $Y_{T}=u(0,X_{T})=\varphi (X_{T})$. That is to say that $%
Y_{t}=u(T-t,X_{t})$ solves the scalar BSDE 
\begin{equation}
dY_{t}=-f_{T}(t,Y_{t},Z_{t})dt+\sum_{j=1}^{d}Z_{t}^{j}dM_{t}^{j}\text{, \ \ }%
Y_{T}=\varphi (X_{T})  \label{bs-2}
\end{equation}
where $f_{T}=f(T-t,y,z)$.

For any fixed $T>0$, let $Y^{T}=\left\{ Y_{t}^{T}:t\in \lbrack 0,T]\right\} $
be the unique solution to the BSDE (\ref{bs-2}) on $(\Omega ,\mathscr{F}%
^{\mu },\mathscr{F}_{t}^{\mu },P^{\mu })$. Since the solution to BSDE is
unique, $Y_{t}^{T}=u(T-t,X_{t})$. In particular $u(T,X_{0})=Y_{0}^{T} $, and
therefore 
\begin{equation}
\int_{R^{d}}u(T,x)\mu (dx)=E^{\mu }\left( Y_{0}^{T}\right) \text{.}
\label{no-l1a}
\end{equation}

The above argument leading to the probabilistic representation (\ref{no-l1a}%
) can not be justified in the case that $a=(a^{ij})$ is only Borel
measurable or the boundary $\partial D$ is only continuous, as in this case, 
$(X_{t})_{t\geq 0}$ is no longer a semimartingale, both (\ref{m-3-a}) and (%
\ref{m-3}) no longer make sense. While, in this case, boundary problem
(strong or weak solutions) to the non-linear PDE (\ref{no-l1}) also need to
be interpreted. On the other hand, the BSDE (\ref{bs-2}), which relies on
only the martingale representation, still make sense, thus the
representation theorems stated in \S3 can be made as the definition of a
solution to (\ref{no-l1}). This is the approach we will carry out.

Consider the initial value problem of the following non-linear parabolic
equation in a bounded domain $D$ with a continuous boundary $\partial D$ 
\begin{equation}
\left( \frac{\partial }{\partial t}-\frac{1}{2}\sum_{i,j=1}^{d}\frac{%
\partial }{\partial x^{j}}a^{ij}(x)\frac{\partial }{\partial x^{i}}\right)
u+f(t,u,\nabla u)=0\text{\ }  \label{n-1}
\end{equation}
subject to the initial and boundary conditions 
\begin{equation*}
u(x,0)=\varphi (x),\text{ \ \ }\left. \frac{\partial }{\partial \nu }%
u(t,\cdot )\right| _{\partial D}=0\text{ for }t>0\text{\ \ }
\end{equation*}
where $a=(a^{ij})$ is Borel measurable, satisfying the uniform ellipticity
condition: 
\begin{equation*}
\lambda \sum_{i=1}^{d}|\xi ^{i}|^{2}\leq \sum_{i,j}^{d}\xi ^{i}\xi
^{j}a^{ij}(x)\leq \lambda ^{-1}\sum_{i=1}^{d}|\xi ^{i}|^{2}\text{ \ \ \ \ \
\ }\forall (\xi ^{i})\in \mathbb{R}^{d},
\end{equation*}
for some constant $\lambda >0$.

\begin{definition}
The functional on $\mathscr{S}_{1}(\overline{D})$ defined by $\mu
\rightarrow E^{\mu }\left\{ Y(t,\mu )_{0}\right\} $, denoted by $u(t,\mu )$,
is called the stochastic solution of \ the initial and boundary problem of (%
\ref{n-1}), where for each $t>0$ and $\mu \in \mathscr{S}_{1}(\overline{D})$%
, $Y(t,\mu )=(Y_{s})_{s\leq t}$ is the unique solution to the BSDE 
\begin{equation}
\left\{ dY_{s}=-f(t-s,Y_{s},Z_{s})ds+\sum_{j=1}^{d}Z_{s}^{j}dM_{s}^{j}\text{%
, \ }Y_{t}=\varphi (X_{t})\right.  \label{bs-3}
\end{equation}
on $(\Omega ,\mathscr{F}^{\mu },\mathscr{F}_{t}^{\mu },X_{t},\theta
_{t},P^{\mu })$.
\end{definition}

As a consequence we have

\begin{theorem}
If $\varphi $ is bounded and Borel measurable on $\overline{D}$, and $f$ is
Lipschitz continuous, then there is a unique stochastic solution to the
non-linear parabolic equation (\ref{n-1})
\end{theorem}

We will study the regularity theory of the stochastic solutions in a
separate paper. On the other hand we would like to derive an alternative
probability representation of the stochastic solution.

Let us apply the approach outlined in \cite{lyons-etc1}. Let $%
Y_{s}=N_{s}-V_{s}$ where 
\begin{equation*}
N_{s}-N_{0}=\sum_{j=1}^{d}\int_{0}^{s}Z_{r}^{j}dM_{r}^{j}\text{.}
\end{equation*}
Then $V=(V_{s})_{s\in \lbrack 0,t]}$ is the unique solution to the
functional differential equation 
\begin{equation}
V_{s}=\int_{0}^{s}f(t-r,Y(V)_{r},Z(V)_{r})dr\text{, \ }V_{0}=0  \label{eq-it}
\end{equation}
where $N(V)_{s}=-E^{\mu }\{\varphi (X_{t})+V_{t}|\mathscr{F}_{s}^{\mu }\}$, 
\begin{equation*}
Y(V)_{s}=E^{\mu }\{\varphi (X_{t})+V_{t}|\mathscr{F}_{s}^{\mu }\}-V_{s}
\end{equation*}
for\ $s\in \lbrack 0,t]$, and $Z(V)$ is given as the density process of $%
N(V) $ in the martingale representation. \ In particular 
\begin{equation*}
Y(V)_{0}=E^{\mu }\{\varphi (X_{t})+V_{t}|\mathscr{F}_{0}^{\mu }\}\text{.}
\end{equation*}
We therefore have the following

\begin{theorem}
Let $\varphi $ be bounded and measurable. For $t>0$, let $V(t)$ be the
unique solution to the functional differential equation (\ref{eq-it}). Then
the stochastic solution to the Neumann boundary problem of the non-linear
PDE (\ref{n-1}) is given by 
\begin{equation}
u(t,\mu )=E^{\mu }\left\{ \varphi (X_{t})+V(t)_{t}\right\} \text{ \ \ \ \ \ }%
\forall \mu \in \mathscr{S}_{1}(\overline{D})\text{.}  \label{a-s1}
\end{equation}
\end{theorem}

\vskip0.4truecm

\noindent {\textbf{Acknowledgements.}} The research of the first author was
supported in part by EPSRC\ grant EP/F029578/1, and by the Oxford-Man
Institute. The second author's research was supported in part by the
National Basic Research Program of China (973 Program) under grant No.
2007CB814904, and a Royal Society Visiting grant.

%\bibliography{myreference}

\noindent {\small \textsc{Z. Qian}}

\noindent{\small Mathematical Institute and Oxford-Man Institute}

\noindent{\small University of Oxford}

\noindent{\small Oxford OX1 3LB, England}

\noindent {\small {Email: \texttt{qianz@maths.ox.ac.uk}}}

\vskip0.3truecm

\noindent {\small \textsc{J. Ying}}

\noindent {\small Institute of Mathematics}

\noindent {\small Fudan University}

\noindent {\small Shanghai, China}

\noindent {\small {Email: \texttt{jgying@fudan.edu.cn}}}

\end{document}